\documentclass[12pt,a4paper]{article}

\usepackage[margin=20mm]{geometry}
\usepackage[T1]{fontenc}
\usepackage[latin2]{inputenc}
\usepackage{ae}
\usepackage{graphicx}
\usepackage[american]{babel}
\usepackage{amsfonts}

\usepackage{amsthm}
\usepackage{amsmath}
\usepackage{amssymb}
\usepackage{wasysym}
\usepackage{paralist}
\usepackage{multirow}
\usepackage{cite}
\usepackage{pifont}
\usepackage{imakeidx}
\makeindex
\makeatletter
\newcommand{\footnoteref}[1]{%
\ltx@ifpackageloaded{hyperref}{
  \ifHy@hyperfootnotes
    \hbox{\hyperref[#1]{%
            \@textsuperscript {\normalfont \ref*{#1}}}}%
  \else
    \hbox{\@textsuperscript {\normalfont \ref*{#1}}}%
  \fi%
}{
    \hbox{\@textsuperscript {\normalfont \ref{#1}}}%
 }%
}
\makeatother

\newtheorem{Theorem}{Theorem}

\newtheorem{Observation}[Theorem]{Observation}

\newcommand{\eps}{\varepsilon}

\usepackage{color}
\definecolor{red}{rgb}{0.8,0,0}
\definecolor{green}{rgb}{0,0.6,0}

\begin{document}

\title{Independent sets and cuts in large-girth regular graphs}
\author{Endre Cs\'oka\thanks{Supported by ERC grants~306493 and 648017.} \\ \normalsize{Alfréd Rényi Institute of Mathematics, Budapest, Hungary} }
\date{}

\maketitle

\begin{abstract}
We present a local algorithm producing an independent set of expected size $0.44533n$ on large-girth 3-regular graphs and $0.40407n$ on large-girth 4-regular graphs. We also construct a cut (or bisection or bipartite subgraph) with $1.34105n$ edges on large-girth 3-regular graphs. These decrease the gaps between the best known upper and lower bounds from $0.0178$ to $0.01$, from $0.0242$ to $0.0123$ and from $0.0724$ to $0.0616$, respectively. We are using local algorithms, therefore, the method also provides upper bounds for the fractional coloring numbers of $1 / 0.44533 \approx 2.24554$ and $1 / 0.40407 \approx 2.4748$ and fractional edge coloring number $1.5 / 1.34105 \approx 1.1185$. Our algorithms are applications of the technique introduced by Hoppen and Wormald.\cite{hoppen2013local}
\end{abstract}

\section{Introduction}

There is a large literature on the relative size of the largest structures with a locally defined property, these questions are also motivated by statistical physics and graph limit theory. The most classical parameter of this kind is the independence ratio, but there is a significant literature on the size of the largest cut and some other parameters, as well.

All of the lower bounds we have are constructive and valid for all large-girth regular graphs. Namely, these largest independent sets and the cuts are constructed by local algorithms, or in other words, by constant-time distributed algorithms. On the other hand, all of the upper bounds are showed for random regular graphs, and these are important open questions whether the random regular graphs have asymptotically the lowest independence ratio, cut ratio, etc. In 2014, Gamarnik and Sudan \cite{gamarnik2014limits} proved that for large $d$, the independence ratio of random $d$-regular graphs cannot be approximated by a local algorithm, not even the $85.4\%$ of it can be achieved. However, the question is open for small $d$, and the best known upper and lower bounds are very close to each other for $d = 3$ and $d = 4$. The question is open for the cut ratio for all $d$.

An independent set in a graph is a subset of the vertices so that no two of them are adjacent. The independence ratio of a graph is the size of the maximum independent set divided by the number of vertices. In 1981, Bollobás\cite{bollobas1981independence} showed that the supremum of the independence ratios of random 3-regular graphs on $n \rightarrow \infty$ vertices is less than $\frac{1}{2}$. In particular, he showed that this is at most $0.4591$. In this paper, Bollobás introduced the configuration model of random graphs, many of the results we will mention are using this tool. In 1982, McKay\cite{mckay} improved the upper bound to $0.45537$, using the same technique but with much more careful calculations. This is still the best upper bound we have.

In 1983, Shearer\cite{shearer1983note} proved the first lower bound on the independence ratio of large-girth 3-regular graphs, this bound was $0.4139$.
In 2008, Hoppen\cite{hoppen2008properties} improved the bound to $0.4328$, his method was improved by Kardo\v{s}, Král' and Volec\cite{kardovs2011fractional}, providing independence ratio $0.4352$.
In 2015 (2013 in arXiv), Csóka, Gerencsér, Harangi and Virág\cite{csoka2015invariant} showed the lower bound of $0.4361$ or $0.438$,\footnote{$\approx 0.4362\pm($proved numeric error bound $0.0001) + ($a nonnegative remaining term, statistical computer simulation showed that this is $\approx 0.0019$, but more than $0.0018)$} using Gaussian wave functions.
Slightly later but with independent techniques (2013 in arXiv), Hoppen and Wormald\cite{hoppen2013local} achieved\footnote{This was the numerical result of a differential equation system. Therefore, there were three best lower bounds: $0.4361$ is rigorously proved, $0.43757$ uses a numerically stable kind of computer computation but with no proved numerical error bound, and $0.438$ uses statistical computer simulation. We consider $0.43757$ to be the most relevant known lower bound prior the present paper. Our result uses the same kind of numerically stable computation.} $0.43757$. They are using the configuration model to build a random graph and an independent set of it at the same time. This is the technique we will use with some additional observations, hereby achieving independence ratio $0.44533$.

For 4-regular random graphs, McKay showed the best upper bound of $0.41635$, and Hoppen and Wormald had the best lower bound of $0.39213$. We improve the lower bound to $0.40407$.

\smallskip

A cut is a 2-coloring of the vertices of a graph, and the size of a cut is the number of edges having endpoints of different colors.
The relative size of the cut is the size of the cut divided by the number of vertices.

The best known upper bound of the relative size of the largest cut for random 3-regular graphs is $1.4026$. It was announced by McKay in 1982 \cite{mckay1982maximum}, its rigorous proof was made by Hladk{\`y} \cite{hladky2006bipartite}. The problem could also be translated to a problem in statistical physics and applying non-rigorous methods suggests that $1.386$ is an upper bound \cite{zdeborova2010conjecture}.

In 1990, Z{\`y}ka proved that the relative size of the maximum cut in cubic graphs with large girth is asymptotically at least $9 / 7 \approx 1.2857$.
In 2002, D{\'\i}az, Do, Serna and Wormald\cite{diaz2002bisection} improved the lower bound to $1.32595$.
In 2012, Kardo{\v{s}}, Kr{\'a}l' and Volec\cite{kardovs2012maximum} constructed $1.33008$.
Then Csóka, Gerencsér, Harangi and Virág introduced the Gauss wave function, and a simple calculation shows that the vertices with positive and negative values make a cut of size $\frac{3}{4} + \frac{3}{\pi} \arctan \frac{1}{\sqrt{2}} \approx 1.33774$. Moreover, this bipartition is not even locally optimal, so it is easy to improve that bound. However, they did not point out this result in their paper, because they were focusing only on the independence ratio.
The technique of Hoppen and Wormald\cite{hoppen2013local} provided a bound of $1.33021$.
With an additional idea, we improve this bound to $1.34105$.

\section{Results on the independence ratio}

First, we show a randomized algorithm that chooses an independent set of expected size at least $0.44531 \cdot n - o(n)$ in a random 3-regular graph on $n$ vertices. Then we will see that this algorithm is local, therefore, this works on all 3-regular graphs with large girth. We will use the technique of Hoppen and Wormald combined with the following observation.

\begin{Observation} \label{obs2}
Assume that a vertex $y$ of a graph $G$ has exactly two neighbors: $x$ and $z$. Let $G'$ be the graph constructed from $G$ by contracting $x$, $y$ and $z$, or in other words, replacing them with a new vertex $v$ with $\deg(v) = \deg(x) + \deg(z) - 2$ and connecting $v$ to the non-deleted neighbors of $x$ and $z$. Then for any independent set $I'$ of $G'$, we can construct an independent set $I$ of $G$ of size $|I| = |I'| + 1$. Namely, if $v \in I'$, then $I = I' \setminus \{v\} \cup \{x, y\}$, otherwise $I = I' \cup \{y\}$.
\end{Observation}

Now we construct a random 3-regular graph and an independent set in it. We use the configuration model introduced by Bollobás\cite{bollobas1981independence}. Namely, in the beginning, we have $n$ vertices of degree 3, but we do not make the pairing of the open edges. Then we construct a random 3-regular graph and an independent set of it at the same time. During the process, we will always have isolated vertices, each of them having some open edges, and only the degree distribution and the number of vertices are changing.

We use the abbreviation $k$\emph{-vertex} for a vertex of degree $k$. If there is a 2-vertex, then we query its neighbours, or in other words, we choose its two neighbors at random with probabilities proportional to the degrees, and we contract these three vertices. If all "survival" vertices (including the new vertices after contraction) have degrees at least 3, then we delete a vertex of the highest degree. This also means that we choose its neighbors at random, and those neighbors lose an open edge. Vertices with degree 0 or 1 will appear in $o(n)$ number of times in expectation, so we can ignore them.
At the end, when only $\eps n$ vertices remain, we choose the empty set as independent set. Then recursively applying Observation~\ref{obs2}, this provides an independent set of a random 3-regular graph of the same size as many contractions we used.

The entire process can be approximated by a quadratic differential equation system on the space of degree distributions. Hoppen and Wormald\cite{hoppen2013local} rigorously proved that this calculation always provides the true asymptotical value, with $o(n)$ error. We cannot solve the system explicitly, but many different versions of calculations provide the same result with 6 digit precision: 0.445312.

We have to add the technical note that in order to apply the theorems of Hoppen and Wormald, instead of using vertices with open edges, we should use
trees with the following bipartite structure. The first class only consists of vertices with degree 2, and both neighbors of each of these vertices go to the second class. The edges leaving from the second class are either go to the first class or these are open edges. (Therefore, the second class is larger then the first class by exactly 1 vertex.) Now the degree of a vertex is replaced with the total number of the open edges of a tree. Deleting a vertex is replaced with adding the first class to the independent set and deleting the tree including the randomly chosen pairs of the open edges. Contracting a 2-vertex is replaced with connecting 3 trees, and hereby all vertices in the middle tree change to the other class.

We give a sketch about what is happening during the process. At the very beginning, we delete 1 (or maybe $\eps n$) number of 3-vertices. Then, there are positive fractions of 3- and 4-vertices only. At each step, we delete a 4-vertex (potentially creating 2-vertices) and we make a contraction with all vertices of degree 2 (potentially creating 5- and 6-vertices), and we delete all vertices of degree 5 or 6 (potentially creating 2-vertices), until none of them exist. Let us consider an arbitrary point of the process when the proportion of 4-vertices is $\lambda$ and $\mu = \frac{4\lambda}{3 + \lambda}$. Now deleting a random open edge changes a 4-vertex to a 3-vertex with probability $\mu$ and a 3-vertex to a 2-vertex with probability $1 - \mu$. A 2-vertex contracts two 3-vertices with probability $(1 - \mu)^2$, a 3- and a 4-vertex to a 5-vertex with probability $2\mu(1 - \mu)$ and two 4-vertex to a 6-vertex with probability $\mu^2$. After a contraction, deleting a these 5- and 6-vertices deletes $5 \cdot 2\mu(1 - \mu) + 6 \mu^2 = (10 - 4\mu)\mu$ open edges, in expectation. Therefore, if $(10 - 4\mu)\mu < 1$, then if we delete an $\eps n$ fraction of the 4-vertices, then it induces a total of $\frac{4}{1 - (10 - 4\mu)\mu} \eps n$ deletions of random open edges, $\frac{4(1 - \mu)}{1 - (10 - 4\mu)\mu} \eps n$ contractions, $\frac{4\mu - 4(1 - \mu)(3 - 2\mu)}{1 - (10 - 4\mu)\mu} \eps n$ change in the number of 3-vertices and $\frac{4(1 - \mu)^3 - 4\mu - 1}{1 - (10 - 4\mu)\mu} \eps n$ change in the number of 4-vertices.

After a point where $(10 - 4\mu)\mu = 1$, we no longer delete 4-vertices. But there are positive fractions of 3-, 4- and 5-vertices, we immediately use contractions at all 2-vertices and delete all 6-, 7- and 8-vertices, and if none of them exist, then we delete a 5-vertex. After two more phases, the process converges to a stationary degree distribution $d_3 \approx 0.55$, $d_4 \approx 0.26$, $d_5 \approx 0.131$, $d_6 \approx 0.055$ and $d_7 \approx 0.004$. This means that if we run the process with this degree distribution, then the number of vertices of each kind decreases proportionally, their ratio remains the same.
We delete approximately $24\%$, $50\%$, $25\%$ and $0.5\%$ of the vertices in the four phases, respectively (in the sense that a contraction counts two deletions). We note that the total proportion of 5-, 6- and 7-vertices, compared to the \emph{original} number of vertices remains below $0.016$ throughout the process.

We made a minor improvement that in the first phase, when we choose a vertex of degree 4, we choose its four neighbors, and if we find two of them with degree 4, then we delete those two vertices and apply contraction at the original vertex. Otherwise we just delete the original vertex. Unfortunately, this improved the resulting independence ratio from $0.445312$ to only $0.445327$. The computer calculation can be found at http://codepad.org/dtBqmn2H .

Now we sketch how this process can be applied in any graph with large girth. This is a classical technique from local algorithms, also described in \cite{hoppen2013local}. We approximate this continuous process by a finite-round discretized version, similarly to what the numeric calculation does. Namely, in each round, we delete $\eps n$ vertices of the highest degrees, and then we make the contractions at the 2-vertices. When we delete some of the vertices of the same degree, then we use independent randomization at each vertex, with the same probability. Say, in the second phase, we delete all 6- 7- and 8-vertices, and we delete each 5-vertex with a small fixed probability. We stop the process when the survival graph has less than $\eps n$ vertices, in expectation. Hereby, the resulting independence ratio can be arbitrarily close to the ratio provided by the continuous process. This way, after a constant number of rounds, the result at each vertex depends only on the randomization in its constant-radius neighborhood. Therefore, the probability of being in the independent set depends only on the constant radius neighborhood of the vertex. Thus, we can construct a random independent set of expected size $0.445327$ in any graph with girth large enough.

Almost the same process provides independence ratio $0.404073$ on 4-regular graphs. The only difference is that if we have no vertex of degree 2 or at least 6, then we do the following. We choose a 3-vertex $v$, we choose its neighbors, and if all of them are 3-vertices, then we delete $v$, otherwise we delete its neighbor with the highest degree, and apply contraction at $v$. The computer calculation can be found at http://codepad.org/lEBUQJu2 .

\section{Results on the size of the maximum cut}

With almost the same techniques, we construct a cut with relative size $1.34105$ on large-girth 3-regular graphs. However, we will temporarily have multiple-vertex components. We will color the vertices with the output colors red and green. Except that we will color some vertices $v$ the "temporary output color" white and choose a neighbor $w$ of $v$ in the survival graph declaring that after the end of the process, if $w$ will be red, then $v$ will be recolored green, and if $w$ will be green, then $v$ will be red. Whenever we assign an output color to the vertex, we immediately remove it from the survival graph. At the very last step, we will not use the color white, this guarantees that all vertices will have the final color red or green. The role of red and green will be symmetric and their distributions will be the same throughout the process.

For the types of the survival vertices, we use the notation $[$list of output colors R, W or G of neighbors$]$. E.g. $[]$ denotes a 3-vertex, $[R]$ denotes a 2-vertex in the survival graph with a red third neighbor, $[GG]$ is a 1-vertex with two other green neighbors. Components are denoted by e.g. $[]-[R]-[G]$, all of our components will be paths, and vertices of type $T$ will be called $T$-vertices.

\begin{Observation} \label{white}
For a $[W]$-vertex $v$, let us delete it and if both of its neighbors are in the component, then we connect these two neighbors and we flip the colors from red to green and vice versa on one side of $v$. E.g. $[R]-[W]-[R]$ is transformed to $[R]-[G]$ or $[G]-[R]$, and $[]-[W]-[R]-[W]-[W]-[R]-[W]-[G]-[W]$ is transformed to $[]-[G]-[G]-[G]$ or $[]-[R]-[R]-[R]$. When a component with flipped colors expands, then this flipping extends to the new vertices. The final color of $v$ will be the opposite of the color of one of its neighbors. Therefore, assuming that the entire process is symmetric to red and green, we can delete all vertices of type $[W]$ with each time just increasing the number of good edges by 1.
\end{Observation}

\begin{Observation} \label{rrr}
A subgraph $[R]-[R]-[R]$ can be replaced with $[R]$ so that the output color of the left and right vertices will be the same as the color of the new vertex, and the color of the middle vertex will be the opposite of it. With this transformation, we should increase the number of good edges (which is in the cut) by 3 and the number of bad edges by 1.
\end{Observation}

The following table shows the list of components in the survival graph and the action we make with them. We apply these rules for any subgraphs of the components of the survival graph. Also, we apply the rules with exchanging red and green.

\begin{center}
\begin{tabular}[center]{|c|c|c|}
\hline
\textbf{Component} & \textbf{After action} & \textbf{Action in words} \\
\hline \hline
$[RR]$ or $[RW]$ & $G$ & color it green \\
\hline
$[RG]$ or $[WW]$ & $W$ & color it white, choose the only survival neighbor \\
\hline
$[W]$ & nothing & apply Observation~\ref{white} \\
\hline
$[G]-[R]$ & $R, G$ & color them red and green \\
\hline
$[]-[R]-[]$ & $[G], G, [G]$ & color the middle vertex green \\
\hline
$[]-[R]-[R]-[]$ & $[G], G, [RG]-[]$ & color one of the middle vertices green \\
\hline
$[R]-[R]-[R]$ & $\cong [R]$ & apply Observation~\ref{rrr} \\
\hline
$[]-[R]$ & $[]-[R]-?$ & query the other neighbor of $[R]$ \\
\hline
$[]-[R]-[R]$ & $[]-[R]-[R]-?$ & query the other neighbor of the second $[R]$ \\
\hline
\hline
$\multirow{2}{*}{[R]}$ & lowest priority: & IF we only have components $[]$ and $[R]$ and $[G]$: \\
 & $[R]-?$ & we query a neighbor of a random $[R]$ or $[G]$ \\
\hline
\end{tabular}
\end{center}

We apply the last step with $[R]$ or $[G]$ with equal probabilities, each time independently. The distribution of the components at an intermediate step is $\lambda$ proportion of $[R]$, $\lambda$ proportion of $[G]$ and $1 - 2\lambda$ proportion of $[]$, ignoring the $[W]$-vertices. All the other components are removed immediately. In order to follow the process, we have to solve the following equation system describing what happens while the number of $[]$-vertices decreases by 1. $p = \frac{2\lambda}{2 + 2\lambda}$, $c$ refers to the expected number of times the different types of actions are used, $v_R$ is the change in the total number of $[R]$- and $[G]$-vertices, $r$ is the number of times we match an open edge of a vertex with output color red or green to a random open edge, $w$ is the same for white, $g$ and $b$ are the number of good and bad edges (in or not in the cut).

\begin{align*}
v_R &= - c_R - 2p \cdot (c_R + r + c_{3R} + c_{3RR} + c_{RR} + w) + (1 - 2p) \cdot (r + c_{3RR}) + (2 - 3p) \cdot c_{3R} + p \cdot c_{RR}
\\ 0 &= - r + p \cdot (2 \cdot c_R + r + c_{3R} + c_{3RR} + c_{RR} + 2 \cdot w)
\\ 0 &= - c_{3R} + (1 - 2p) \cdot c_R + p \cdot c_{3RR}
\\ 0 &= - c_{3RR} + p \cdot c_{3R} + (1 - 2p) \cdot c_{RR}
\\ 0 &= - c_{RR} + p \cdot c_R
\\ 0 &= - w + p \cdot (r + c_{RR})
\\ - 1 &= - (1 - 2p) \cdot (c_R + r + c_{3R} + c_{3RR} + c_{RR} + w)
\\ g &= 3p \cdot c_R + 4p \cdot r + (1 + p) \cdot c_{3R} + (4 + p) \cdot c_{3RR} + 8p \cdot c_{RR} + w
\\ b &= p \cdot r + c_{3RR} + 2p \cdot c_{RR}
\end{align*}

The results multiplied by $2 - 4 p - 4 p^2 + 8 p^3 + 2 p^4 - 4 p^5$ are $v_R = 1 - 8 p + 4 p^2 + 8 p^3 + 3 p^4 - 10 p^5$, $g = 1 + 8 p - 11 p^2 - 6 p^3 + 12 p^5$ and $b = p (1 - p)^2 (2 + p + 2 p^2)$. From here, http://codepad.org/f66tIj1o shows the computer calculation.

\bibliography{reference-independence}

\newpage

\section*{Appendix}

We show here the program codes we used. In these calculations, the proportions which should be 0 are fluctuating between $\pm \eps$. We note that many different ways of calculations were tested, and all of them agree with 6 digits precision. These codes run at codepad.org, in C++.

\subsection*{Independence ratio for large-girth 3-regular graphs}

http://codepad.org/dtBqmn2H

\begin{verbatim}
int main()
{
  const long double eps = 0.0000000063;
  cout << setprecision(9);
  const int deg = 8; // Maximum degree + 1
  long double v[deg];
  for (int i = 0; i < deg; ++i)
    v[i] = 0.;
  v[3] = 1;
  long double independent = 0;
  long double erase = 0;
  int round = 0;
  while (v[3] > eps)
  {
    ++round;
    long double s = 0;
    for (int i = 3; i < deg; ++i)
      if (v[i] > eps)
        s += i * v[i];
    if (erase > eps)
    {
      long double r = (erase + eps) / s;
      for (int i = 3; i < deg; ++i)
        if (v[i] > eps)
        {
          long double del = r * i * v[i];
          v[i] -= del;
          v[i - 1] += del;
        }
      erase = -eps;
    }
    s = 0;
    for (int i = 3; i < deg; ++i)
      if (v[i] > eps)
        s += i * v[i];
    if (v[2] > eps)
    {
      long double r = (v[2] + eps) / s;
      long double add[2 * deg - 3];
      for (int i = 3; i < 2 * deg - 3; ++i)
        add[i] = 0;
      for (int i = 3; i < deg; ++i)
        for (int j = 3; j < deg; ++j)
          add[i + j - 2] += i * v[i] * j * v[j];
      independent += v[2] + eps;
      v[2] = -eps;
      for (int i = 3; i < deg; ++i)
        v[i] += r * (add[i] / s - 2 * i * v[i]);
      for (int i = deg; i < 2 * deg - 3; ++i)
        erase += i * r * add[i] / s;
    }
    int max = deg - 1;
    while ((max > 4) && (v[max] < eps))
      --max;

    v[max] -= 2 * eps;
    erase += 2 * max * eps;

    if (max == 4)
    {
      long double p4444 = 2 * eps * pow(4 * v[4] / s, 4);
      v[3] -= 4 * p4444;
      erase += 12 * p4444;
      independent += p4444;

      long double p4443 = 8 * eps * pow(4 * v[4] / s, 3) * 3 * v[3] / s;
      v[3] -= 3 * p4443;
      v[2] -= p4443;
      erase += 11 * p4443;
      independent += p4443;

      long double p4433 = 12 * eps * pow(12 * v[4] * v[3] / s / s, 2);
      v[4] += p4433;
      v[3] -= 2 * p4433;
      v[2] -= 2 * p4433;
      erase += 6 * p4433;
      independent += p4433;
    }
    if ((round % 131072 == 0) || (v[3] <= eps))
    {
      cout << "independent = " << independent;
      for (int i = 3; i < deg; ++i)
        cout <<", v[" << i << "] = " << v[i];
      cout << endl;
    }
  }
  return 0;
}
\end{verbatim}

\newpage

\subsection*{Independence ratio for large-girth 4-regular graphs}

http://codepad.org/lEBUQJu2

\begin{verbatim}
int main()
{
  const long double eps = 0.00000001;
  const int deg = 8; // Maximum degree + 1
  cout << setprecision(9);
  long double v[deg];
  for (int i = 0; i < deg; ++i)
    v[i] = 0.;
  v[4] = 1;
  long double independent = 0;
  long double erase = 0;
  int round = 0;
  while (v[4] > eps)
  {
    ++round;
    long double s = 0;
    for (int i = 3; i < deg; ++i)
      if (v[i] > eps)
        s += i * v[i];
    if (erase > eps)
    {
      long double r = (erase + eps) / s;
      for (int i = 3; i < deg; ++i)
        if (v[i] > eps)
        {
          long double del = r * i * v[i];
          v[i] -= del;
          v[i - 1] += del;
        }
      erase = -eps;
    }
    s = 0;
    for (int i = 3; i < deg; ++i)
      if (v[i] > eps)
        s += i * v[i];
    if (v[2] > eps)
    {
      long double r = (v[2] + eps) / s;
      long double add[2 * deg - 3];
      for (int i = 3; i < 2 * deg - 3; ++i)
        add[i] = 0;
      for (int i = 3; i < deg; ++i)
        for (int j = 3; j < deg; ++j)
          add[i + j - 2] += i * v[i] * j * v[j];
      independent += v[2] + eps;
      v[2] = - eps;
      for (int i = 3; i < deg; ++i)
        v[i] += r * (add[i] / s - 2 * i * v[i]);
      for (int i = deg; i < 2 * deg - 3; ++i)
        erase += i * r * add[i] / s;
    }
    int max = deg - 1;
    while ((max > 5) && (v[max] < eps))
      --max;
    if (max > 5)
    {
      v[max] -= 2 * eps;
      erase += 2 * max * eps;
    }
    else
    {
      const long double rat3 = 3 * v[3] / (3 * v[3] + 4 * v[4] + 5 * v[5]);
      const long double rat4 = 4 * v[4] / (3 * v[3] + 4 * v[4] + 5 * v[5]);
      const long double rat5 = 5 * v[5] / (3 * v[3] + 4 * v[4] + 5 * v[5]);
      v[2] += eps * 3 * rat3 * rat3 * rat3;
      v[3] += eps * (-1 - 3 * rat3);
      v[4] += eps * 3 * (- rat4 + rat3 * rat3 * (1 - rat3));
      v[5] += eps * 3 * (- rat5 + rat3 * rat4 * (rat4 + 2 * rat5));
      independent += eps * (1 - rat3 * rat3 * rat3);
      erase += eps * (6 - 12 * rat3 * rat3 + 6 * rat3 * rat3 * rat3 
        + (15 * rat3 * rat4 + 3) * (rat4 + 2 * rat5));
    }
    if ((round % 131072 == 0) || (v[4] <= eps))
    {
      cout << "independent = " << independent;
      for (int i = 3; i < deg; ++i)
        cout <<", v[" << i << "] = " << v[i];
      cout << endl;
    }
  }
  return 0;
}
\end{verbatim}

\newpage

\subsection*{Size of the maximum cut for large-girth 3-regular graphs}

http://codepad.org/f66tIj1o

\begin{verbatim}
int main()
{
  const long double eps = 0.000000011;
  cout << setprecision(9);
  long double good = 0;
  long double bad = 0;
  long double rat3 = 1;
  long double rat2 = 0;
  int round = 0;
  while (rat2 + rat3 > eps)
  {
    ++round;
    const long double q = rat2 / (2 * rat2 + 3 * rat3);
    const long double q2 = pow(q, 2);
    const long double q3 = pow(q, 3);
    const long double q4 = pow(q, 4);
    const long double q5 = pow(q, 5);
    rat2 += eps * (1 - 8 * q + 4 * q2 + 8 * q3 + 3 * q4 - 10 * q5);
    rat3 += eps * (-2 + 4 * q + 4 * q2 - 8 * q3 - 2 * q4 + 4 * q5);
    good += eps * (1 + 8 * q - 11 * q2 - 6 * q3 + 12 * q5);
    bad += eps * q * (1 - q) * (1 - q) * (2 + q + 2 * q2);
    if ((round % 1048576 == 0) || (rat2 + rat3 <= eps))
      cout << "good = " << good << ", bad = " << bad << ", rat3 = " << rat3
        << ", rat2 = " << rat2 << endl;
  }
  return 0;
}
\end{verbatim}

\end{document}